\documentclass[12pt,fleqn]{article}

\usepackage{amsmath}
\usepackage{amsthm}
\usepackage{amssymb}
\usepackage{amsfonts}
\usepackage{epsf}
\usepackage{graphicx}
\usepackage{epstopdf}
\usepackage{hyperref}
\usepackage{color,soul}

\newtheorem{thm}{Theorem}

\theoremstyle{remark}
\newtheorem{rmk}{Remark}

\theoremstyle{definition}

 \DeclareMathOperator\re{{Re}}
 \numberwithin{equation}{section}

\newcommand{\D}{\displaystyle}

\numberwithin{equation}{section}

\newcounter{comment}
\def\a{\alpha}
\def\b{\beta}

\renewcommand{\l}{\lambda}

\begin{document}
\title{Asymptotics of orthogonal polynomials and the Painlev\'e transcendents}
\date{\today}
\author{Dan Dai \\ \\ Department of Mathematics \\ City University of
Hong Kong \\ Hong Kong \\
Email: \texttt{dandai@cityu.edu.hk}}

\maketitle

\begin{center}
  \emph{Dedicated to Professor Mourad Ismail on the occasion of his 70th birthday.}
\end{center}

\vspace{2cm}

\begin{abstract}
  In this survey, we review asymptotic results of some orthogonal polynomials which are relate to the Painlev\'e transcendents. They are obtained by applying Deift-Zhou's nonlinear steepest descent method for Riemann-Hilbert problems. In the last part of this article, we list several open problems.     
\end{abstract}


\vspace{2cm}



\newpage


\section{Introduction}

\subsection{Orthogonal polynomials}

Let $\mu$ be a positive Borel measure on $\mathbb{R}$ and assume that the moments $  \int x^n d\mu(x)$ exist for all $n=0,1,2, \cdots$. Then, there exists a unique sequence of monic orthogonal polynomials $\{\pi_n(x)\}_{n=0}^\infty$,
\begin{equation}
  \pi_n(x) = x^n + \cdots,
\end{equation}
such that
\begin{equation} \label{conti-ortho}
  \int \pi_m(x) \pi_n(x) d\mu(x) = h_n \delta_{m,n}.
\end{equation}
When $d\mu(x)$ = $w(x)dx$ for some continuous function $w(x)$ on an interval, the corresponding polynomials are called \emph{continuous orthogonal polynomials} and the function $w(x)$ is referred to as the weight function. If $\mu$ is a discrete measure, the corresponding polynomials are called \emph{discrete orthogonal polynomials}. For example, when $\mu$ is supported on the integers $k\in \mathbb{Z}$ with masses $w_k$, then instead of  \eqref{conti-ortho}, the polynomials satisfy the following orthogonality relation
\begin{equation} \label{dis-ortho}
  \sum_{k\in \mathbb{Z}} \pi_m(k) \pi_n(k) w_k= h_n \delta_{m,n}.
\end{equation}

It is well-known that orthogonal polynomials play a remarkable role in many areas of mathematics and physics, such as continued fractions, operator theory (Jacobi operators), approximation theory, numerical analysis, quadrature, combinatorics, random matrices, Radon transform, computer tomography, etc. Of course, the list presented here is far from complete.

Besides their important applications, orthogonal polynomials also satisfy a lot of fascinating properties. Let us mention a few of them below. The classical ones, namely \emph{the Jacobi, Laguerre and Hermite polynomials}, satisfy a second-order differential equation of the following form
\begin{equation} \label{intro-ode}
c_2(x) y''(x) + c_1(x) y'(x) + c_0(x) y(x) = 0, 
\end{equation}
where $c_0(x)$, $c_1(x)$ and $c_2(x)$ are polynomials of degree at
most 0, 1 and 2, respectively. From the orthogonality property
(\ref{conti-ortho}), it can be verified that all orthogonal polynomials satisfy a three-term recurrence
relation
\begin{equation} \label{intro-rec}
x \pi_n(x) = \pi_{n+1}(x) + a_n \pi_n(x) + b_n \pi_{n-1}(x), \quad n = 1,
2, \cdots,
\end{equation}
with
\begin{equation}
  \pi_0(x) = 1, \qquad \pi_1(x) = x - a_0,
\end{equation}
where $a_n$ is real and $b_n > 0$ for $n>0$. In addition, the continuous orthogonal polynomials satisfy the well-known Rodrigues formula
\begin{equation}
\pi_n(x) = \frac{1}{c_n w(x)} \frac{d^n}{dx^n} (w(x) [g(x)]^n ),
\end{equation}
where $c_n$ is a constant and $g(x)$ is a polynomial independent
of $n$. Both continuous and discrete orthogonal polynomials have
generating functions
\begin{equation}
F(x,z) = \sum_{n=0}^{\infty} d_n \pi_n(x) z^n.
\end{equation}
Here, for the same polynomials $\pi_n(x)$, the function $F(x,z)$ may
be different by choosing different coefficients $d_n$. From the
generating function and Cauchy's integral formula, one can easily
get an integral representation of the form
\begin{equation} \label{intro-int}
\pi_n(x) = \frac{1}{2 \pi i \, d_n} \int_C \frac{F(x,z)}{z^{n+1}} dz,
\end{equation}
where $C$ is a closed contour surrounding the origin in the positive direction. For more properties of orthogonal polynomials, we
refer to Andrews, Askey and Roy \cite{And:Ask:RoyBook}, Chihara \cite{Chihara:book}, Gautschi \cite{Gautschi:book}, Ismail \cite{Ismbook}, Koekoek, Lesky and Swarttouw \cite{Koek:book}, Levin and Lubinsky \cite{Levin:Lub}, Szeg\H{o} \cite{Szego}, etc.


\subsection{Asymptotic methods}

In the study of orthogonal polynomials, there is a lot of interest in their asymptotic behaviors as the polynomial degree $n$ is large. In the literature, there are several
different methods to derive their asymptotics. For the classical orthogonal polynomials, since they satisfy differential equations of the form (\ref{intro-ode}), one can apply the powerful
asymptotic methods in the differential equation theory, such as the Liouville-Green approximation (also called the WKB approximation) and the turning point theory. For more
details about the differential equation method, we refer to the definitive work of Olver \cite{Olver:book}. However, since not all orthogonal polynomials satisfy differential equations,
this approach becomes inapplicable and people turn to the
integral approach. From the generating function and after
a suitable re-scaling of the variable $x$, one often transforms
the integral (\ref{intro-int}) into the form
\begin{equation}
\frac{1}{2 \pi i} \int_C g(x,z) e^{n f(x,z)} dz.
\end{equation}
Then, based on the properties of the phase function $f(x,z)$ and its
saddle points, one may apply the method of steepest descents to derive
the asymptotic expansion of the orthogonal polynomials. For more
information about the integral approach, we refer to Wong \cite{Wong:intBook}.

Besides the differential equation method and the integral approach mentioned above, researchers have made significant progress on two novel methods in the past twenty years, namely the difference equation method and the Riemann-Hilbert(RH) approach. Starting from the three-term recurrence relation (\ref{intro-rec}), Wong and his colleagues developed a turning point theory for second-order difference equations in a series of papers \cite{Cao:Li,Wang:Wong2003,Wang:Wong2005}. Their ideas and asymptotic results are similar to those in the differential equation theory. But the asymptotic analysis is far more complicated than the differential cases due to the discrete feature. And the turning point
theory for difference equations is completely non-trivial. Nowadays, the difference equation theory is viewed as an analogue of the classical asymptotic theory for linear second-order differential equations; see the survey article by Wong \cite{Wong:AA2014}.

\subsection{The Riemann-Hilbert approach} \label{sec:op-rhp}

The RH approach is based on the key observation that orthogonal polynomials are connected with a RH problem as follows (cf. Fokas, Its and Kitaev \cite{fik1992}).

\medskip

\noindent\textbf{The RH problem for $Y(z)$:}

\medskip

Assume $\mu$ is supported on the real line and $d\mu(x)$ = $w(x)dx$ in \eqref{conti-ortho}. The problem is to determine a $2 \times 2$ matrix-valued function $Y: \mathbb{C} \to \mathbb{C}^{2 \times 2}$ such that the following holds.

\begin{description}
  \item(Y1) $Y(z)$ is analytic in
  $\mathbb{C}\backslash \mathbb{R}$;

  \item(Y2) Let $ Y_+(x)$ and $Y_-(x)$ denote the limiting values of $Y (z)$ as $z$ approaches $x\in\mathbb{R}$ from
the upper and lower half plane, respectively. Then, $Y(z)$  satisfies the jump condition
  \begin{equation}\label{Y-jump}
  Y_+(x)=Y_-(x) \left(
                               \begin{array}{cc}
                                 1 & w(x) \\
                                 0 & 1 \\
                                 \end{array}
                             \right),
\qquad x\in \mathbb{R};
\end{equation}

 \item(Y3) The asymptotic behavior of $Y(z)$  at infinity is
  \begin{equation}\label{Y-infty}Y(z)=\left (I+O\left (  1 /z\right )\right )\left(
                               \begin{array}{cc}
                                 z^n & 0 \\
                                 0 & z^{-n} \\
                               \end{array}
                             \right),\quad \mbox{as}\quad z\rightarrow
                             \infty .
  \end{equation}

\end{description}

By virtue of the Plemelj formula and Liouville's theorem, one has the following results.
\begin{thm} \label{thm-fokas} (Fokas, Its and Kitaev \cite{fik1992})
  The unique solution to the above RH problem is given by
  \begin{equation}\label{Y-solution}
Y(z)= \left (\begin{array}{cc}
\pi_n(z)& \D\frac 1 {2\pi i}
\int_{\mathbb{R}}\frac {\pi_n(s) w(s) }{s-z} ds\\[0.35cm]
-2\pi i \gamma_{n-1}^2 \;\pi_{n-1}(z)& -  \D  \gamma_{n-1}^2\;
\int_{\mathbb{R}}\frac {\pi_{n-1}(s) w(s) }{s-z} ds \end{array} \right ),
\end{equation}
where  $\pi_n(z)$ is the monic orthogonal polynomial with respect to the weight
$w(x)$, and $\gamma_n$ is the leading coefficient of the orthonormal polynomial.

\end{thm}

\begin{rmk}
  If the weight function $w(x)$ is supported on certain interval $(a,b)$ instead of $\mathbb{R}$, similar results also hold. The only difference is that one needs to add extra conditions at the endpoints  $a,b$ to ensure the uniqueness of the RH problem. For example, see Kuijlaars et  al. \cite{Kui:McL:Ass:Van2004} where $(a,b) = (-1,1)$.
\end{rmk}

\begin{rmk} \label{rmk:discrete}
  One can also formulate a RH problem for the discrete orthogonal polynomials. In this case, the above RH problem becomes an  interpolation one and  the jump condition (Y2) is replaced by residue conditions; see Baik et al. \cite{Baik:Kri:McL:Miller}.
\end{rmk}

On the basis of Theorem \ref{thm-fokas}, Deift and Zhou et al. \cite{dkmv1,dkmv2} developed a very powerful nonlinear steepest descent method to derive the asymptotics of $\pi_n(z)$ from the RH problem for $Y$; see also Deift \cite{deift}. The idea is to obtain, via a series of invertible transformations $Y \to T \to S \to R $, the RH problem for $R$ whose jump is close to the identity matrix.
\begin{itemize}
    \item $Y\to T$ is to rescale the variable $z$, and normalize large-$z$ condition (Y3). As a result, $T(z)$ solves a RH problem with oscillatory jumps, normalized at infinity.

    \item $T\to S$ is to deform the contour. In the meantime, the oscillatory jumps are factorized and $S(z)$ solves a RH problem without oscillation.

    \item $S\to R$, the final transformation, includes global parametrix construction as well as local parametrix constructions near critical points (usually they are also the self-intersection points of the contour). Then, $R$ satisfies a RH problem with all jumps close to $I$  for large polynomial degree $n$.
\end{itemize}
The asymptotic expansion for $R$ can be derived on the whole complex plane. Tracing back, the uniform asymptotics of the orthogonal polynomials is obtained as $n \to \infty$. Since the establishment of Deift-Zhou nonlinear steepest descent method, the RH approach has been successfully applied to study asymptotics of various orthogonal polynomials, as well as their application in proving universality results in random matrix theory, see \cite{Baik:Dei:Joh,Baik:Kri:McL:Miller,Claeys:Kuij:Van,Dei:Its:Kra2011,dkmv1,Kui:McL:Ass:Van2004,vanlessen} for a very incomplete list from the literature.

It is well-known that, the asymptotic expansions of classical orthogonal polynomials usually involve Airy functions, Bessel functions and trigonometric functions, see \cite[Sec. 18.15]{nist}. All of these functions satisfy linear ordinary differential equations(ODEs). In the past a few years, people realize that, when a parameter tends to certain critical values in weight functions, some new functions will appear in the asymptotic expansions. These functions are related to the Painlev\'e equations, which are nonlinear second-order ODEs. The RH approach plays a crucial role in deriving these Painlev\'e asymptotics. Indeed, such kind of results have not been obtained by any other methods so far. In this article, we will review asymptotic results of some orthogonal polynomials which are related to the Painlev\'e equations. However, due to the space limitation, it is impossible for us to cover all of the results. For example, we don't mention asymptotics of multiple orthogonal polynomials, which have very important applications in random matrix theory. As a consequence, if any related results are not covered in this survey, the omission certainly does not reflect on the importance of the omitted works.

The rest of the article is arranged as follow. In Sec. \ref{sec:Pain}, we first briefly introduce the definitions of the Painlev\'e equations as well as some of their properties. Then, we focus on the RH problems associated with the Painlev\'e equations. In Sec. \ref{sec:Asy}, after providing readers with some ideas about why and how the Painlev\'e asymptotics will appear, we review the asymptotic results of orthogonal polynomials which are related to PI-PV in the literature. In the last section, we list several problems in this research area.


\section{Painlev\'e equations} \label{sec:Pain}

In the literature, people are interested in nonlinear ordinary differential equations of the following form
\begin{equation} \label{orig-ode}
  w_{xx} = F(x, w, w_x),
\end{equation}
where $F$ is a function meromorphic in $x$ and rational in $w(x)$
and $w'(x)$. Usually the solutions of these equations have movable
singularities, which means that the positions of possible
singularities depend on the initial conditions of the equations. At
the turn of the twentieth century, Painlev\'e, Gambier, Fuchs et al. (1893--1906) classified the equations \eqref{orig-ode} whose solutions are free from movable branch points and
essential singularities. This property is called \emph{the Painlev\'e property} now. It turns out that the
equations satisfying the Painlev\'e property can be reduced to six canonical forms, which are now known as the Painlev\'e differential equations listed below.
\begin{eqnarray*}
    \hbox{PI: \quad} w_{xx} & = & 6 w^2 + x \\
    \hbox{PII: \quad} w_{xx} & = & 2 w^3  + x w  - \alpha \\
    \hbox{PIII: \quad} w_{xx} & = & \frac{1}{w} w_x^2 - \frac{1}{x} w_x  + \frac{1}{x}
    (\alpha w^2 + \beta) + \gamma w^3 + \frac{\delta}{w}
\end{eqnarray*}
\begin{eqnarray*}
    \hbox{PIV: \quad} w_{xx} & = & \frac{1}{2w} w_x^2 + \frac{3}{2} w^3 +
    4 x w^2 + 2(x^2 - \alpha)w + \frac{\beta}{w}, \\
    \hbox{PV: \quad} w_{xx} & = & \left( \frac{1}{2w} + \frac{1}{w-1} \right) w_x^2 -\frac{1}{x}w_x + \frac{(w-1)^2}{x^2} \left(\a w + \frac{\b}{w}\right) \\
    && +\frac{\gamma w}{x} + \frac{\delta w (w+1)}{w-1}, \\
    \hbox{PVI: \quad} w_{xx} & = & \frac{1}{2} \left( \frac{1}{w} + \frac{1}{w-1} + \frac{1}{w-x} \right) w_x^2  - \left( \frac{1}{x}  + \frac{1}{x-1} + \frac{1}{w-x} \right) w_x \\
     &&+ \frac{w(w-1)(w-x)}{x^2(x-1)^2} \left( \a + \frac{\b x}{w^2} + \frac{\gamma (x-1)}{(w-1)^2} + \frac{\delta x (x-1)}{(w-x)^2} \right),
\end{eqnarray*}
where $\alpha, \beta, \gamma$ and $\delta$ are constants. The solutions of PI-PVI are called the \emph{Painlev\'e transcendents}.

One can see that the study of the Painlev\'e equations originates from a purely mathematical point of view. However, it turns out that the Painlev\'e equations satisfy a lot of nice properties and have applications in wide areas of physics. Let us briefly mention several properties first. For example, the Painlev\'e equations appear as similarity reductions of some integrable nonlinear PDEs. For some special parameters, they possess exact solutions which are given in terms of rational functions, algebraic functions or classical special functions. For PII-PVI, the Painlev\'e transcendents satisfy \emph{B\"{a}cklund transforms}, which can be viewed as nonlinear recurrence relations. Moreover, similar to the role that the classical special functions play in linear physics, the Painlev\'e transcendents appear in many areas of nonlinear physics, such as statistical mechanics, random matrix theory, quantum gravity and quantum field theory, nonlinear optics and fibre optics, etc. Nowadays, the Painlev\'e transcendents are viewed as nonlinear analogues of the classical special functions and included in the latest version of mathematical handbook \cite[Chap. 32]{nist}; also see Clarkson \cite{Clarkson:Lec2006} and references therein for more information.

\subsection{Lax pair and Riemann-Hilbert problems} \label{sec:pain-rhp}

Although the Painlev\'e equations are nonlinear, each of them can be
expressed as the compatibility condition of a system of
\emph{linear} differential equations (Lax pair) in the following
form
\begin{equation} \label{laxpair}
    \frac{\partial \Psi}{\partial \l} = A(\l, x) \Psi, \quad  \frac{\partial \Psi}{\partial x} = B(\l, x) \Psi.
\end{equation}
Here, $\Psi( \l,x)$, $A(\l,x)$ and $B(\l,x)$ are $2 \times 2$ matrix
functions, $A(\l,x)$ and $B(\l,x)$ are rational functions in $\l$.
The compatibility condition is $\D\frac{\partial^2 \Psi}{\partial \l \partial x}
= \frac{\partial^2 \Psi}{\partial x \partial \l}$, which means
\begin{equation}
    A_x - B_\l + AB - BA = 0.
\end{equation}
The Lax pair (\ref{laxpair}) is greatly helpful in studying properties of the Painlev\'e transcendents. One of the direct outcomes is that, from the ordinary differential equation theory
and the rational property of $A(\l,x)$, it is possible for us to formulate a RH problem for the sectionally holomorphic function
$\Psi( \l,x)$ in the complex $\l$-plane, with $x$ appearing as a parameter. The solution of this
RH problem is associated with some special solutions of the Painlev\'{e} equations if the jump matrices are specified; see Fokas et al. \cite{fikn} for comprehensive information about this method.

Let us take PI as a concrete example to illustrate the relation between a RH problem and the Painlev\'{e} transcendents. In \cite{Kapaev}, Kapaev formulates a model  RH  problem for $\Psi(\l)= \Psi(\l; x)$ as follows.

\begin{itemize}

  \item[(a)] $\Psi(\l; x)$ is analytic for $\l \in \mathbb{C} \setminus \Gamma_\Psi$, where
  \begin{equation}
    \Gamma_\Psi= \gamma_{-2} \cup \gamma_{-1} \cup \gamma_{1} \cup \gamma_{2} \cup \gamma^*
  \end{equation}
  is illustrated in Figure \ref{contour-p1-model}.

  \begin{figure}[h]
 \begin{center}
   \includegraphics[width=7.5cm]{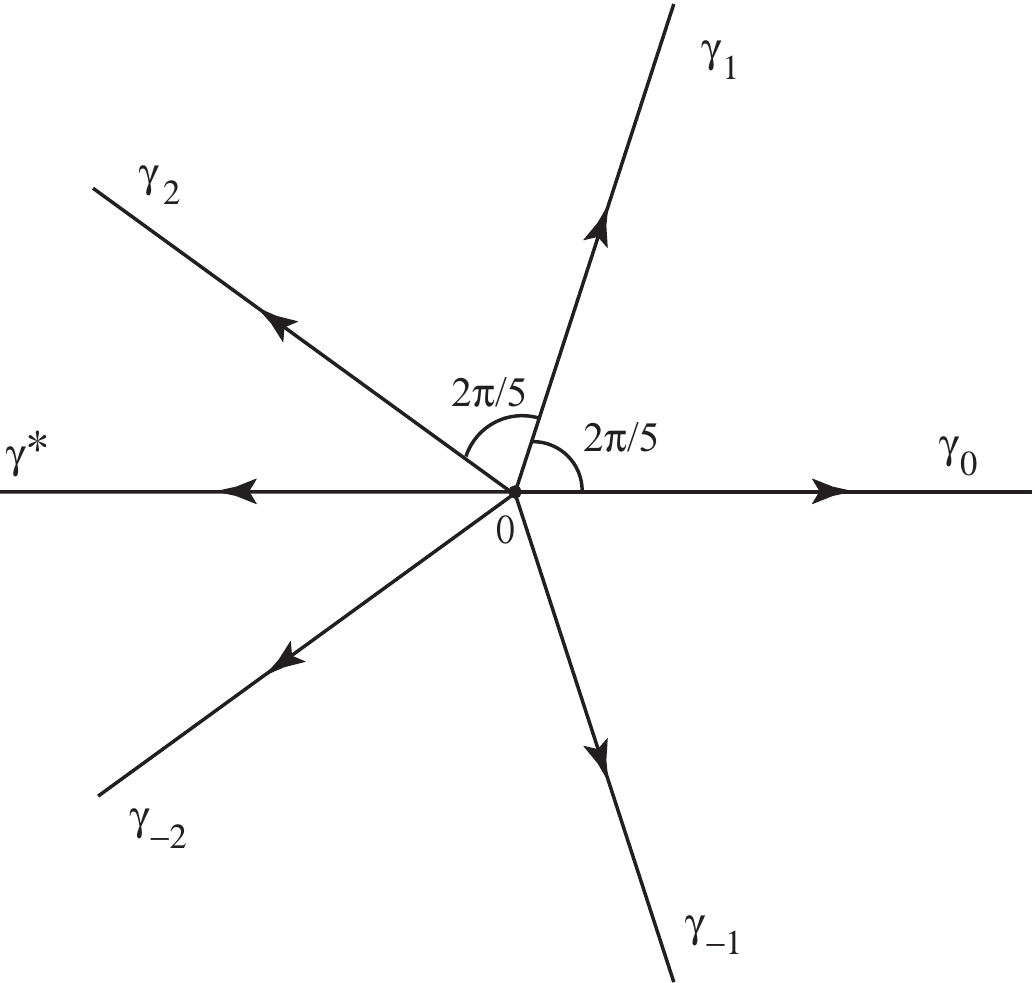} \end{center}
  \caption{The contour $\Gamma_\Psi$ associated with the Painlev\'e I equation}
 \label{contour-p1-model}
\end{figure}

  \item[(b)] Let $\Psi_{\pm} (\l; x)$ denote the limiting values of $\Psi(\l; s)$ as $\l$ tends to the contour $\Gamma_\Psi$ from the left and right sides, respectively. Then, $\Psi(\l; s)$ satisfies the following jump conditions
   \begin{equation}\label{Psi--jump}
  \Psi_+(\l; x)= \Psi_-(\l; x)
 \left\{\begin{array}{ll}
          \left(
                               \begin{array}{cc}
                                 1 & s_k \\
                               0& 1 \\
                                 \end{array}
                             \right), & z \in \gamma_k, \ k=\pm1; \\[.4cm]
           \left(
                               \begin{array}{cc}
                                 1 & 0 \\
                                 s_k& 1 \\
                                 \end{array}
                             \right), & z \in \gamma_k, \ k=\pm 2;  \\[.4cm]
          \left(
                               \begin{array}{cc}
                                0 &-i\\
                               -i&0 \\
                                 \end{array}
                             \right),  &  z\in \gamma^*,
        \end{array}\right .
 \end{equation}
        where the complex constants $s_k$'s are the so-called \emph{Stokes multipliers}. They satisfy the following relation
        \begin{equation} \label{stoke-multilier}
          1 + s_k s_{k+1} =-i s_{k+3} , \qquad  s_{k+5} = s_k , \qquad  k \in \mathbb{Z}.
        \end{equation}

   \item[(c)] As $\l \to \infty$, $\Psi(\l; x)$ satisfies the asymptotic condition
  \begin{equation}\label{Psi-infty}
    \Psi(\l; x)= \l^{\frac{1}{4} \sigma_3}  \frac{\sigma_3+\sigma_1}{\sqrt{2}}\left(I + \frac{\Psi_{-1}(x)}{\sqrt{\l}} + \frac{\Psi_{-2}(x)}{\l} + O(\l^{-\frac 32})\right)  e^{\theta(\l,x) \sigma_3}
    \end{equation}
    for $\arg \l\in (-\pi, \pi)$, where
    \begin{equation} \label{theta-def}
      \theta(\l,x)=\frac 45 \l^{\frac 52}+x\l^{\frac 12},
    \end{equation}
    $\sigma_1$ and $\sigma_3$ are the Pauli matrices
    \begin{equation*}
      \sigma_1 = \begin{pmatrix}
        0 & 1 \\ 1 & 0
      \end{pmatrix}, \qquad \sigma_3 = \begin{pmatrix}
        1 & 0 \\ 0 & -1
      \end{pmatrix}.
    \end{equation*}

\end{itemize}
From the above RH problem, one can prove that
\begin{equation} \label{tronquee-sol}
  w(x) = 2(\Psi_{-2} (x))_{12}
\end{equation}
is a solution of PI, where $\Psi_{-2} (x)$ is given in \eqref{Psi-infty}. Moreover, it is known that, for every set of Stokes multipliers $s_k$ satisfying \eqref{stoke-multilier}, there exists a unique solution of PI and vice versa. In particular, by choosing
\begin{equation}
  s_0 = 0, \quad s_1=(1-\alpha)i, \quad s_{-1}=\alpha i, \quad s_{\pm2} = i
\end{equation}
with $\alpha$ being a complex constant, we get a special solution $w_\alpha(x)$ of PI.
One can even derive the asymptotic expansions of $w_\alpha(z)$ from the above RH problem when the independent variable $z$ is complex. More precisely, Kapaev \cite{Kapaev} showed that $w_\alpha(z)$ is the so-called \emph{tronqu\'ee} solution of PI whose asymptotic behavior is given by
\begin{equation} \label{p1-a}
w_\alpha(z)=w_0(z)+\frac {\alpha i}{\sqrt{\pi}}2^{-\frac {11}8}(-3z)^{-\frac 18}\exp\left[-\frac 152^{\frac {11}4}3^{\frac 14}(-z)^{\frac 54}\right]\left (1+O(z^{-\frac 38})\right )
\end{equation}
as $z\to \infty$ and $\arg z =\arg(-z)+\pi    \in[\frac 35\pi,\pi]$. Here $w_0(z)$ is the \emph{tritronqu\'ee} solution satisfying
\begin{equation}\label{Painleve I}
w_0(z)\sim\sqrt{-z/6}\left[ 1+ \sum_{k=1}^{\infty}a_k(-z)^{-5k/2} \right] \quad \textrm{as } z \to \infty, \ -\frac{\pi}{5}<\arg z< \frac{7\pi}{5},
\end{equation}
where the coefficients $a_k$ can be determined recursively. In general, the PI transcendents are meromorphic functions and possess infinitely many poles in the complex plane. The tronqu\'ee and tritronqu\'ee solutions mentioned above satisfy the property that the solution $w_\alpha(z)$ is pole-free for $z$ in some sectors of the complex $z$-plane. For more information about the tritronqu\'ee solutions of PI, we refer to Costin et al. \cite{Cos:Huang:Tan}, Joshi and Kitaev \cite{Jos:Kit}, as well as references therein.

\begin{rmk} \label{rmk:p6}
  The RH problems for PII-PV can be found in Fokas et al. \cite[Chap. 5]{fikn}. About PVI, its Lax pair possesses four regular singular points, which are more than the numbers of singular points of any other Painlev\'e equations; see Jimbo and Miwa \cite[Appendix C]{Jim:Miw}. As a consequence, the RH problem for PVI is a little more complicated. For example, from the associated Lax pair, Its, Lisovyy and Prokhorov constructed a RH problem for PVI in  \cite[Sec. 3]{Its:Lis:Pro2016}. Then, they use this problem to study asymptotic behaviors of the corresponding tau functions.
\end{rmk}


\section{The Painlev\'e asymptotics} \label{sec:Asy}

From Section \ref{sec:op-rhp} and Section \ref{sec:pain-rhp}, one can see that both orthogonal polynomials and the Painlev\'e transcendents satisfy $2 \times 2$ RH problems. Then, the natural question is that whether there exist any relations between these two types of RH problems. The answer is positive. Indeed, this is also the reason why the Painlev\'e type asymptotics appear.

Let us recall the Deift-Zhou nonlinear steepest descent method. The final transform $S \to R$ involves parametrix constructions, which are essentially constructions of approximate solutions to the RH problem $S$ when the parameter $n$ is large. Under certain special situations, the RH problems for orthogonal polynomials near some critical points are very similar to those for the Painlev\'e equations. Therefore, the RH problems for the Painlev\'e equations (for example the one in Section \ref{sec:pain-rhp}) is adopted to construct the local parametrix in the transform $S \to R$. As a consequence, the Painlev\'e functions appear in the final asymptotic expansions of orthogonal polynomials.

To make a more clear explanation about how the Painlev\'e type asymptotics appear, we need the limiting zero distribution of orthogonal polynomials. Given a well-behaved weight function, this distribution can be obtained explicitly from the potential theory; see Saff and Totik \cite{Saff:Totik}. Let us rescale the variable and put the weight function $w(x)$ in the form of $e^{-nV(x)}$. The function $V(x)$ is also called \emph{potential}. For simplicity, we first assume the weight function is supported on the whole real axis.
It is well-known that the limiting zero distribution is the so-called \emph{equilibrium measure} $\nu^*$, which is the unique minimizer of following energy functional among all Borel probability measures $\nu$ on $\mathbb{R}$
\begin{equation} \label{energy-minimize}
  E_V(\nu):=\int\int \log \frac{1}{|x-y|} d\nu(x) d\nu(y) + \int V(x) d\nu(x).
\end{equation}
If $V(x)$ is real analytic, then
\begin{equation}
  \textrm{supp}\, \nu^* = \bigcup_{i=1}^k [a_i, b_i], \qquad a_1<b_1 < \cdots < a_k < b_k,
\end{equation}
where $k$ is a finite number. Moreover, the density of $\nu^*$ is given in the following form
\begin{equation} \label{e-measure1}
  \rho(x) = \frac{d\nu^*}{dx} = \begin{cases}
  h(x) \sqrt{(x-a_1)(b_1-x) }, & \textrm{if } k =1 \\
    h(x) \sqrt{(x-a_1)(x-b_1) \cdots (x-a_k)(b_k-x)}, & \textrm{if } k \geq 2
  \end{cases}
\end{equation}
for $x \in \textrm{supp}\, \nu^*$,
where $h(x)$ is a real analytic function and strictly positive for $x \in \textrm{supp}\, \nu^*$; see Deift et al. \cite{Dei:Kri:McL1998}. Note that the density function $\rho(x)$ is strictly positive in the support and vanishes like a square root at the endpoints. Under this situation, the corresponding asymptotic expansions of orthogonal polynomials are given in terms of trigonometric functions in the compact subset of $\textrm{supp}\, \nu^*$ and Airy functions in the neighbourhood of $a_i$ and $b_i$, respectively. If the weight function is supported not on the whole real axis but on an interval (for example $[a,\infty)$), the density function $\rho(x)$ may blow up with an exponent $-1/2$ near the bounded endpoint $a$, namely,
\begin{equation} \label{e-measure2}
  \rho(x) \sim \frac{1}{\sqrt{x-a}} \qquad \textrm{as } x \to a+.
\end{equation}
Then, the asymptotic expansions of orthogonal polynomials involve Bessel functions near the endpoint $a$.

The equilibrium measure depends on the weight function of orthogonal polynomials. For some special weights, the corresponding equilibrium measures are no longer in the form of \eqref{e-measure1} or \eqref{e-measure2}. Consequently, some new functions, such as the Painlev\'e functions, will emerge in the asymptotic expansions of orthogonal polynomials.

\subsection{PI asymptotics} \label{sec:p1}

When the density of an associated equilibrium measure vanishes with an exponent 3/2 at an endpoint of its support, the local parametrix construction will involve RH problems of PI.
Note that this type of vanishing is impossible in the case of usual orthogonality with respect to exponential weights on the real line. One needs to consider the  \emph{non-Hermitian orthogonality} which we will explain below with a concrete example.

In \cite{Dui:Kuij}, Duits and Kuijlaars consider a varying quartic weight
\begin{equation} \label{duits-weight}
  w_N(x):= e^{-NV_t(x)} = e^{-N(tx^4/4 + x^2/2)} \qquad \textrm{for } t<0.
\end{equation}
Obviously, when $t<0$, all of the moments
\begin{equation} \label{quart-moments}
  \mu_{k,N} = \int_\mathbb{R} x^k e^{-N(tx^4/4 + x^2/2)} dx
\end{equation}
do not exist. As a consequence, the corresponding orthogonal polynomials are not well-defined as well. To ensure the convergence of the above integral, one needs to change the original integration contour $\mathbb{R}$ to a contour $\Gamma$ in the complex plane, such that $\re V_t(z) \to +\infty$ as $z \to \infty$ along the contour $\Gamma.$ By choosing $\Gamma $ to be $\{re^{\pi i /4}, r\in\mathbb{R}\}$, 
we turn to study the polynomials satisfying the following non-Hermitian orthogonality
\begin{equation}
  \int_{\Gamma} \pi_{n,N}(z) z^k e^{-N(tz^4/4 + z^2/2)} dz = 0, \qquad \textrm{for } k =0, 1, \cdots, n-1.
\end{equation}
It is interesting to note that the first RH problem for orthogonal polynomials formulated by  Fokas, Its and Kitaev \cite{fik1992} is related to the non-Hermitian orthogonality instead of the usual one. However, because the bilinear form
\begin{equation}
  <p,q>=\int_\Gamma p(z) q(z) e^{-N(tz^4/4 + z^2/2)} dz
\end{equation}
is not positive definite, the polynomials $\pi_{n,N}(z)$ may not exist for some $n$. One of the results in \cite{Dui:Kuij} is that they proved the existence of $\pi_{n,n}(z)$ and $\pi_{n\pm 1,n}(z)$ when $n$ is large enough.

When $t<0$, it still makes sense to consider the equilibrium problem in \eqref{energy-minimize}. Indeed, similar to the case $t \geq 0$, the density of the equilibrium measure can also calculated explicitly:
\begin{equation}
  \rho(x) = \frac{t}{2 \pi} (x^2 - d_t^2) \sqrt{ c_t^2 - x^2} \qquad \textrm{for } x \in [-c_t, c_t],
\end{equation}
where $t_{cr} < t <0$ with $t_{cr}= -1/12$. Here $c_t$ and $d_t$ are two constants depending on $t$ and $0 <c_t< d_t$. Note that the above density function is similar to that in \eqref{e-measure1}. Therefore, like what we have achieved in the classical cases, trigonometric and Airy type asymptotics will also appear. New phenomenon occurs as $t \to t_{cr}+$, where both $c_t$ and $d_t$ tend to $\sqrt{8}$ . Then, the density becomes
\begin{equation} \label{duits-em}
  \rho_{cr}(x) = \frac{1}{24 \pi} (8-x^2)^{3/2}, \qquad x \in [-\sqrt{8}, \sqrt{8}],
\end{equation}
which vanishes with an exponent 3/2 at endpoints. Taking a \emph{double scaling limit} as $n \to \infty$ and $t \to t_{cr}$, the PI asymptotics appear. To give readers some ideas about the Painlev\'e type asymptotics, we quote Duits and Kuijlaars' main results below. However, the details of other results will not be included in this article due to the space limitation.

\begin{thm} \label{thm-duits} (Duits and Kuijlaars \cite[Thm. 1.2]{Dui:Kuij})
  Let $t$ vary with $n$ such that
  \begin{equation}
    n^{4/5}(t + 1/12) = -c_1x, \qquad c_1 = 2^{-9/5}3^{-6/5},
  \end{equation}
  remains fixed, where $x$ is not a pole of $w_1(x)$ in \eqref{p1-a}. Then, for large enough $n$, the recurrence coefficients $b_{n,n}(t)$ associated with the orthogonal polynomials exist and satisfy
  \begin{equation}
    b_{n,n}(t) = 2 - 2c_2 w_1(x) n^{-2/5} + O(n^{-3/5}), \qquad  c_2 = 2^{3/5}3^{2/5},
  \end{equation}
  as $n\to\infty$. The above expansion holds uniformly for $x$ in compact subsets of $\mathbb{R}$ not containing any of the poles of $w_1(x)$.
\end{thm}

\begin{rmk}
  To simplify the formulas and make them easier to understand, we choose the parameters $\a =\b = 1$ in Theorem 1.2 of \cite{Dui:Kuij}. In this case, the other recurrence coefficient $a_{n,n}(t)$ in \eqref{intro-rec} is equal to 0 for all $n$.
\end{rmk}

It is well-known that the Painlev\'e transcendents possess infinitely many poles in the complex plane in general. In the above theorem, to ensure the validity of the results, the variable $x$ is required to be bounded away from the poles of the tronqu\'ee $w_\a(x)$ in \eqref{p1-a}. Recently, through a more delicate \emph{triple scaling limit}, Bertola and Tovbis \cite{Ber:Tov} successfully obtain the asymptotics near the poles of $w_\a(x)$.

Similar situation like that in \eqref{duits-em} also appears for the cubic potentials
\begin{equation}
  V_t(z) = tz^3 + z^2/2 \qquad \textrm{for } t<0.
\end{equation}In \cite{Ble:Dea1,Ble:Dea}, Bleher and Dea\~{n}o derived the PI asymptotics of the corresponding orthogonal polynomials as well as their applications in random matrix theory. Recently, inspired by the study of Wigner time-delay in mathematical physics (cf. \cite{texier,Wigner1955}), Xu, Dai and Zhao \cite{xdz2016} considered a singular potential as follows
\begin{equation}
  V_t(z) = z- \log z+ t/z \qquad \textrm{for } t<0.
\end{equation}
Note that the above potential has a pole at $z=0$, which is different from the polynomial potentials considered before. Applying the Deift-Zhou steepest descent analysis, we also obtained the PI asymptotics for the corresponding orthogonal polynomials.

\begin{rmk}
  PI asymptotics appear only when the non-Hermitian orthogonality is considered. The reason is that, for the usual orthogonality, the density of the equilibrium measure can only vanish at an endpoint with an exponent $(4k + 1)/2 $ with $k \in \mathbb{N} \cup \{0\}$; see \cite{Dei:Kri:McL1998JAT,Kui:McL2000} for more details. The generic case $k = 0$ gives the Airy-type asymptotics. The first critical case $k = 1$ yields asymptotic expansion in terms of the second member of the Painlev\'e I hierarchy; see Claeys and Vanlessen \cite{Claeys:Van2007}. For $k \geq 2$, it is expected that the asymptotic expansions will involve higher members of the Painlev\'e I hierarchy.
\end{rmk}

\subsection{PII asymptotics}

The first PII asymptotics for orthogonal polynomials appear in the paper \cite{Baik:Dei:Joh} by Baik, Deift and Johansson, where they studied the length of the longest increasing subsequence of a random permutation of $\{1,2,\cdots ,n \}$. In their paper, they need to know the asymptotic behaviour of certain \emph{orthogonal polynomials on the unit circle}(OPUC). However, the following quartic potential may be a simpler example to understand. We will come back to the PII asymptotics about OPUC at the end of this section.

\begin{figure}[h]
 \begin{center}
\includegraphics[width=430pt]{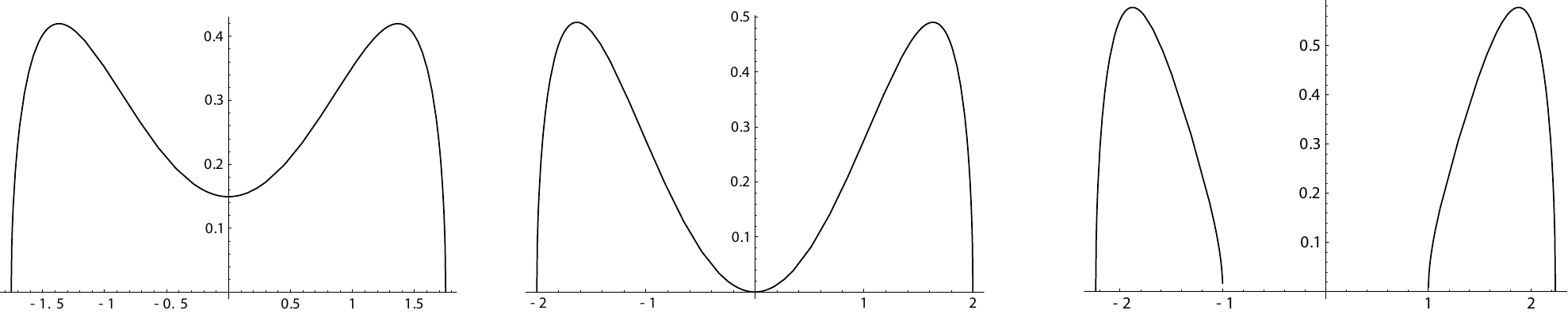}
\end{center}
\caption{From the left to right, the density $\rho(x)$ when $t>-2$, $t=-2$ and $t<-2$, respectively.} \label{fig-quart}
\end{figure}

In \cite{Ble:Its1}, Bleher and Its considered the following varying quartic weight
\begin{equation} \label{ble-its-weight}
  w_N(x):= e^{-NV_t(x)} = e^{-N(x^4/4 + tx^2/2)} , \qquad x \in \mathbb{R}.
\end{equation}
Since the leading coefficient of $V_t(x)$ is positive, all of the moments $\mu_{k,N}$ given in \eqref{quart-moments} exist and corresponding orthogonal polynomials satisfy the usual orthogonality. One may compare the difference between the two quartic weights in \eqref{duits-weight} and \eqref{ble-its-weight}. For the quartic potential in \eqref{ble-its-weight}, the corresponding equilibrium measure can also be computed explicitly. 
Depending on the parameter $t$ in \eqref{ble-its-weight}, we have
\begin{equation} \label{quart-em1}
  \rho(x) = \frac{x^2  + b_t}{2 \pi} \sqrt{a_t^2-x^2}, \qquad x \in [-a_t, a_t], \qquad \textrm{when} \quad t > -2,
\end{equation}
with $a_t,b_t > 0$;
\begin{equation} \label{quart-em-cr}
  \rho_{cr}(x) = \frac{x^2}{2 \pi} \sqrt{4-x^2}, \qquad x \in [-2, 2], \qquad \textrm{when} \quad t = t_{cr} = -2,
\end{equation}
and
\begin{equation} \label{quart-em2}
  \rho(x) = \frac{|x|}{2 \pi} \sqrt{(a_t^2 - x^2)(x^2 - b_t^2)}, \qquad 0 < b_t \leq |x| \leq a_t, \qquad \textrm{when} \quad t < -2,
\end{equation}
with $a_t = \sqrt{2-t}$ and $b_t = \sqrt{-2-t}$; see Figure \ref{fig-quart}.

When $t>-2$ and $t<-2$, the density functions in \eqref{quart-em1} and \eqref{quart-em2} are similar to that in \eqref{e-measure1}. But when $t=-2$, the density function $\rho_{cr}(x)$ in \eqref{quart-em-cr} vanishes quadratically at an interior point of its support. Moreover, the case $t=-2$ indicates a transition where the support of $\rho(x)$ splits from one cut (when $t>-2$) to two cuts (when $t<-2$). Consequently, in the neighbourhood of 0, a PII RH problem is needed to construct the local parametrix and the PII asymptotics follow. Note that the special PII transcendent involved here is the famous Hastings-McLeod solution for homogeneous PII equation(the constant term $\a = 0$ in PII); see \cite{Ble:Its1}. For more properties about the Hastings-McLeod solution, as well as other special solutions of homogeneous PII, we refer to Deift and Zhou \cite{Dei:Zhou1995} and references therein.

Soon, the work of Bleher and Its \cite{Ble:Its1} was generalized by Claeys and Kuijlaars \cite{claeys:Kuijlaars} and  Claeys, Kuijlaars and Vanlessen \cite{Claeys:Kuij:Van}. Instead of the quartic potential in \eqref{ble-its-weight}, it was shown in \cite{claeys:Kuijlaars,Claeys:Kuij:Van} that similar results also hold for more general potential $V(x)$, as long as the corresponding equilibrium measure satisfies the similar properties illustrated in Figure \ref{fig-quart}. In addition, if an algebraic singularity is added at the point where that cut-split takes place (for example, one may include an extra factor $|x|^\alpha$ in the quartic weight \eqref{ble-its-weight}), then we will obtain asymptotics related to the inhomogeneous PII equation; see \cite{claeys:Kuijlaar2008,Claeys:Kuij:Van}.

On the other hand, when studying polynomials orthogonal with respect to the perturbed Hermite weight, people found out that their asymptotics are related to the Painlev\'e XXXIV equation; see Its, Kuijlaars and \"{O}stensson \cite{its:kui:ost} for an algebraic perturbation and Xu and Zhao \cite{Xu:Zhao2011} for a jump perturbation. Note that Painlev\'e XXXIV can be reduced to PII through a simple transform.

Finally, let us go back to OPUC mentioned at the beginning of this section. Let $p_n(z) = \kappa_n z^n + \cdots $, $\kappa_n > 0$, be a polynomial  orthonormal
with respect to the weight $w(z)$ on the unit circle, they satisfy the following orthogonality relation
\begin{equation}
  \frac{1}{2\pi i} \int_{|z| = 1} p_n(z)\, \overline{p_m(z)} \, w(z) \frac{dz}{z} = \delta_{m,n}.
\end{equation}
Sometimes, the above formula is rewritten as
\begin{equation}
  \frac{1}{2\pi } \int_{- \pi}^\pi  p_n(e^{i\theta}) \, \overline{p_m(e^{i\theta})} \, w(e^{i\theta}) d\theta = \delta_{m,n}.
\end{equation}
Instead of \eqref{intro-rec} in the real case, OPUCs satisfy a recurrence relation as follows
\begin{equation}
  \kappa_n p_{n+1}(z) = \kappa_{n+1} z \, p_{n}(z) + p_{n+1}(0) p_{n}^*(z),
\end{equation}
where $p_n^*(z) = z^n \, \overline{p_n(1/\bar{z})}$. From the orthogonality condition, it is also possible to formulate a RH problem for OPUC as follows.
\begin{description}
  \item(Y1) $Y(z)$ is analytic in   $\mathbb{C}\backslash \Sigma$, where $\Sigma$ is  the unit circle oriented counterclockwise.

  \item(Y2) $Y(z)$  satisfies the jump condition
  \begin{equation}
  Y_+(z)=Y_-(z) \left(
                               \begin{array}{cc}
                                 1 & w(z)/z^n \\
                                 0 & 1 \\
                                 \end{array}
                             \right),
\qquad z\in \Sigma.
\end{equation}

 \item(Y3) The asymptotic behavior of $Y(z)$  at infinity is
  \begin{equation}
  Y(z)=\left (I+O\left (  1 /z\right )\right )\left(
                               \begin{array}{cc}
                                 z^n & 0 \\
                                 0 & z^{-n} \\
                               \end{array}
                             \right),\quad \mbox{as}\quad z\rightarrow
                             \infty .
  \end{equation}
\end{description}
For more information and properties about OPUC, we refer to Simon's treatise \cite{Simon:book1,Simon:book2}.

For the OPUC appearing in Baik, Deift and Johansson \cite{Baik:Dei:Joh}, they are orthogonal about the varying weight
\begin{equation}
  w_N(z):= e^{-NV_t(z)} = e^{\frac{tN}{2} (z + z^{-1})}, \qquad z \in \Sigma.
\end{equation}
The associated equilibrium measure is supported on the whole circle or part of it, depending on the parameter $t$. When $t$ achieves the critical value  $t_{cr}=1$, the measure becomes
\begin{equation}
  d\mu_{cr}(\theta) = \frac{1}{2\pi} (1 + \cos \theta) d\theta, \qquad \theta\in [0, 2 \pi].
\end{equation}
One can see that the above density vanishes quadratically at $\theta= \pi$, which is similar to the real quartic case \eqref{quart-em-cr}. Based on the above properties of the equilibrium measure, PII asymptotics appear under a double scaling limit as $n \to \infty$ and $t \to t_{cr}$. It should be emphasized that, this is the first time that the Painlev\'e RH problems emerge in the local parametrix construction of the RH analysis. Since the work of Baik, Deift and Johansson \cite{Baik:Dei:Joh}, more and more Painlev\'e asymptotics are discovered while studying orthogonal polynomials with different weight functions.

\subsection{PIII asymptotics}

Regarding PIII, because both 0 and $\infty$ are irregular singular points in the first equation of its Lax pair \eqref{laxpair}, then the corresponding RH problem possesses an \emph{essential singularity} at 0. While perturbing the classical Hermite and Laguerre weight with some essential singularities,
\begin{eqnarray}
    w(x;t) &=& e^{ -\frac{t}{2x^2} - \frac{x^2}{2}}, \qquad t > 0,  \ x \in \mathbb{R}, \label{Mo-weight} \\
    w(x;t)&= &x^{\alpha}e^{-x - \frac{t}{x}},~~~x\in (0, \infty),~~t>0,~~\alpha>0, \label{xdz-weight}
\end{eqnarray}
Brightmore, Mezzadri and Mo \cite{bmm} and Xu, Dai and Zhao \cite{xdz2014,xdz2015} obtained the PIII asymptotics for polynomials orthogonal with respect to the above weight functions, respectively. Later, Atkin, Claeys and Mezzadri \cite{At:Cla:Mez} generalized \eqref{xdz-weight} and studied the weight
\begin{equation}
  w(x;t) = x^\alpha e^{-x -(\frac{t}{x})^k}, \qquad k \in \mathbb{N}.
\end{equation}
In the end, the asymptotic expansions in \cite{At:Cla:Mez} are given in terms of solutions to a hierarchy of the PIII equations. On the other hand, Xu and Zhao \cite{Xu:Zhao2015} and Zeng, Xu and Zhao \cite{Zeng:Xu:Zhao2015} considered a modified Jacobi weight of the form
\begin{equation}
   w(x;t) = (1 - x^ 2 )^\b (t^2 - x^2 )^\a , \qquad x \in (-1,1)
\end{equation}
with $\b > -1$, $\a + \b > -1$ and $t>1$. Their asymptotics is essentially related to the generalized PV equation. After a M\"{o}bius transformation, this equation can be transformed to a PIII.

\subsection{PIV asymptotics}

To obtain the PIV asymptotics, like Sec. \ref{sec:p1} for the PI cases, we need to change the orthogonality interval on the real line to a contour in the complex plane and consider non-Hermitian orthogonality. In \cite{Dai:Kui}, Dai and Kuijlaars considered the following weight
\begin{equation}
  z^{-N + \nu} e^{-Nz} (z-1)^{2b}.
\end{equation}
The equilibrium measure is supported on the so-called Szeg\H{o} curve $\mathcal{S}$
\begin{equation}
  \mathcal{S} := \{ z \in \mathbb{C}:  |z e^{1-z} | = 1 \quad \textrm{and} \quad |z| \leq 1 \}.
\end{equation}
Moreover, the density of the equilibrium measure vanishes linearly at the point $z=1$. Note that such kind of vanishing behavior is impossible for the usual orthogonality. For the usual case, if the density vanishes at an interior point $z_0$ of its support, it must behaves like $|z-z_0|^{2k}$ for $k\in \mathbb{N}$. When $k=1$, we have a density similar to that in \eqref{quart-em-cr} and obtain the PII asymptotics. For the current case, we achieved the PIV asymptotics instead.

\subsection{PV asymptotics}

Note that the original PV equation possesses three singularities: 0, 1 and $\infty$. The corresponding RH problem also includes three singularities, which are more than those in the RH problems for PI-PIV. Consequently, the contours and the jump conditions in the RH problem for PV become more complicated. On the other hand, if the asymptotics of orthogonal polynomials involve PV, their weight functions are supposed to have some singularities. Moreover, these singularities should influence each other when the parameters in the weight function vary. In \cite{Claeys:Its:Kra}, Claeys, Its and Krasovsky considered OPUC whose weight function possesses a Fisher-Hartwig singularity, which combines a jump-type and a root-type singularity. Taking a double scaling limits, they derived the PV asymptotics for the related OPUC. See also Claeys and Krasovsky \cite{Claeys:Kra2015} for
two merging Fisher-Hartwig singularities for OPUC. Similar situations also exist for polynomials orthogonal on real intervals. For example, Claeys and Fahs \cite{Claeys:Fahs} studied weight functions with merging algebraic singularities; Xu and Zhao \cite{Xu:Zhao2013,Xu:Zhao2015} considered modified Jacobi weights where the algebraic singularity tends to the endpoint of the orthogonal interval. For all these cases, the PV asymptotics come into play.


\section{Further problems}

In this section, we list some unsolved problems in this area.

\medskip

\noindent\textbf{1. The PVI asymptotics}

\medskip

To the best of our knowledge, the PVI asymptotics for orthogonal polynomials have never appeared in the literature. Using the isomonodromy method introduced in Fokas et al. \cite{fikn}, one can derive RH problems for PVI from its Lax pair; see Remark \ref{rmk:p6}. Note that, the RH problems for Painlev\'e equations may not be unique. If one finds a nice RH problem for PVI,  this problem can be employed to construct local parametrices for orthogonal polynomials with certain special weight functions. Then,  the PVI asymptotics will be obtained accordingly.


\medskip

\noindent\textbf{2. Discrete orthogonal polynomials}

\medskip

Although these polynomials satisfy a discrete orthogonality in \eqref{dis-ortho}, it is still possible to formulate a discrete RH problem (also called interpolation problem) for them; see Remark \ref{rmk:discrete}. Starting from this problem, Baik et al. \cite{Baik:Kri:McL:Miller} developed Deift and Zhou's method to derive the asymptotic expansions for discrete orthogonal polynomials. Note that these polynomials possess some special features. For example, there exist the so-called \emph{saturated regions}, where the density of the limiting zero distribution reaches the maximum (the density of the orthogonality nodes distribution) in some intervals. This leads to an \emph{upper constraint} for the equilibrium measure in the energy minimization problem \eqref{energy-minimize}. However, despite these differences, similar asymptotic expansions are also obtained, which involve trigonometric functions and the Airy functions; for example, see Baik et al. \cite{Baik:Kri:McL:Miller}. Recently, while studying nonintersecting Brownian motions on the half-line, Liechty \cite{Liechty2012} considered discrete orthogonal polynomials with respect to the following weight
\begin{equation}
  w(x_k) = \exp\left(-\frac{ n\pi^2 t }{2} x^2 \right ), \qquad x_k = \frac{k-\alpha}{n}, \textrm{ for } k \in \mathbb{Z} \ \textrm{ and } \  \a \in [-1/2, 1/2].
\end{equation}
Note that there exists a critical value $t_{cr}=1$ when the upper constraint of the equilibrium measure is just active. In a double scaling limit as $n \to \infty$ and $t \to 1$, the PII asymptotics for the orthogonal polynomials were derived in \cite{Liechty2012}. See also a subsequent work by Liechty and Wang \cite{Lie:Wang2016} where a class of similar discrete orthogonal polynomials is studied. Now the question is: will there appear any other Painlev\'e asymptotics for discrete orthogonal polynomials?

\medskip

\noindent\textbf{3. $q$-orthogonal polynomials}

\medskip

In Ismail's monograph \cite{Ismbook}, there is a large part related to $q$-series and $q$-orthogonal polynomials. With his colleagues, Ismail has also studied their asymptotics, which usually involve the $q$-Airy functions and Jacobi theta functions; for example see Ismail \cite{Ism-IMRN2005}, Ismail and Li \cite{Ism:Li2013}, Ismail and Zhang \cite{Ism:Zhang2007}, etc.. It is well-known that both $q$-orthogonal polynomials and $q$-Airy functions are $q$-generalizations of the corresponding classical ones. In the literature, there is also a considerable amount of work about $q$-Painlev\'e equations. It would be very interesting if some connections between $q$-orthogonal polynomials and
$q$-Painlev\'e equations can be established.


\section*{Acknowledgements}

The author is grateful for Professors Shuai-Xia Xu and Yu-Qiu Zhao for their
helpful comments and discussions.

The present work was partially supported by grants from the Research Grants Council of the
Hong Kong Special Administrative Region, China (Project No. CityU 11300814, CityU 11300115).


\end{document}